\documentclass[12pt]{amsart}

\usepackage{amsmath,amssymb,amsthm,mathtools}
\usepackage[T1]{fontenc}
\usepackage{lmodern,microtype}
\usepackage[hidelinks]{hyperref}

\newtheorem{theorem}{Theorem}
\newtheorem{lemma}[theorem]{Lemma}
\newtheorem{proposition}[theorem]{Proposition}
\newtheorem*{corollary}{Corollary}

\newcommand{\D}{\mathbb D}
\newcommand{\C}{\mathbb C}
\newcommand{\OO}{\mathcal O}
\newcommand{\LL}{\mathcal L}
\newcommand{\PSH}{\operatorname{PSH}}
\newcommand{\dist}{\operatorname{dist}}
\newcommand{\re}{\operatorname{Re}}

\title[Lempert and Green functions near non-semipositive points]
{Lempert and pluricomplex Green functions near non-semipositive boundary points}

\author{Armen Edigarian}
\address{A. Edigarian\\Department of Mathematical Analysis\\Institute of Mathematics\\
Faculty of Mathematics and Computer Science\\Jagiellonian University\\
ul. Prof. S. \L{}ojasiewicza 6, 30-348 Krak\'ow, Poland}
\email{armen.edigarian@uj.edu.pl}
\thanks{A. Edigarian was funded by the National Science Centre, Poland, under the Weave UNISONO grant no.~UMO-2025/07/Y/ST1/00146.}

\author{Nikolai Nikolov}
\address{N. Nikolov\\Institute of Mathematics and Informatics\\Bulgarian Academy
of Sciences\\Acad. G. Bonchev 8, 1113 Sofia, Bulgaria
\vspace{1mm}
\newline Faculty of Information Sciences\\
State University of Library Studies and Information Technologies\\
Shipchenski prohod 69A, 1574 Sofia, Bulgaria}
\email{nik@math.bas.bg}
\thanks{N. Nikolov was partially supported by the Bulgarian National Science Fund,
Ministry of Education and Science of Bulgaria under contract KP-06-N82/6.}

\date{}

\begin{document}

\begin{abstract}
We study the two-step Lempert function and the pluricomplex Green function near a
non-semipositive boundary point of a bounded domain $D\subset\mathbb{C}^n$. We prove a
sharp upper estimate for the two-step Lempert function in terms of the local boundary
scale $H_D$, and obtain matching estimates for the pluricomplex Green function. As a
consequence, these quantities have the same boundary behaviour as the Kobayashi
distance. For the usual Lempert function, we prove an upper estimate of order
$\sqrt{H_D}$, which implies the optimal local $1/4$-H\"older estimate, and
determine its precise boundary behaviour along the real normal.
\end{abstract}
\maketitle

\section{Introduction and main results}

For a domain $D\subset\C^n$ and a real-valued function
$h:D\times D\to\mathbb R$, we write $h^*:=\tanh h$ and
\[
 h^{(2)}(z,w):=\inf_{u\in D}\bigl(h(z,u)+h(u,w)\bigr),
 \qquad z,w\in D.
\]
For background on the functions used below, see, e.g., \cite{JP}.
Define the Lempert function by
\[
 l_D(z,w)
 :=\operatorname{arctanh}\inf\{\,|\alpha|:\exists f\in\OO(\D,D),\
 f(0)=z,\ f(\alpha)=w\,\}.
\]
Recall that the Kobayashi pseudodistance $k_D$ is the chain pseudodistance
generated by $l_D$; equivalently, it is the largest pseudodistance on $D$ which
does not exceed $l_D$. Thus $l_D^{(2)}$ is the two-step version of this chain
construction. We also use
the Carath\'eodory pseudodistance
\[
 c_D(z,w):=\operatorname{arctanh}\sup\{\,|f(w)|:f\in\OO(D,\D),\ f(z)=0\,\}.
\]
It is the smallest holomorphically contractible pseudodistance with the usual
normalization on $\D$. Hence
\begin{equation}\label{eq:chain}
 c_D^*\le k_D^*
 \le (l_D^{(2)})^*
 \le (l_D^*)^{(2)}
 \le l_D^*.
\end{equation}

The additive pluricomplex Green function with pole at $z$ is
\[
 G_D(z,w):=\sup\{u(w):u\in\PSH(D),\ u<0,\
 u(\zeta)\le\log|\zeta-z|+O(1)\text{ as }\zeta\to z\},
\]
and we use its bounded version
\[
 g_D:=\exp G_D.
\]
We use the standard fact that $G_D(z,\cdot)$ is plurisubharmonic (with the
possible value $-\infty$); see, e.g., \cite{EdiGreen,JP}.
Every $f\in\OO(D,\D)$ with $f(z)=0$ gives an admissible competitor
$\log|f|$ in the definition of $G_D(z,\cdot)$, while holomorphic
contractibility and $g_{\D}(0,\alpha)=|\alpha|$ give the opposite bound. Thus
\begin{equation}\label{eq:green-lempert}
 c_D^*\le g_D\le l_D^*.
\end{equation}
In general there is no two-sided comparison between $g_D$ and $k_D^*$.
The pluricomplex Green function need not be symmetric.

The purpose of the results below is to describe the behaviour of these functions
near non-semipositive boundary points. When both points approach the same such point,
the pluricomplex Green function, the Kobayashi distance, and the two-step Lempert
quantities have a common local scale.  The Lempert function has a different
single-disc behaviour: we obtain the upper estimate $l_D^*\lesssim\sqrt{H_D}$
and hence the optimal local H\"older estimate of exponent $1/4$, as well as
qualitative uniform bounds and a sharp normal estimate.

Let $p\in\partial D$ be a $C^{1,1}$-smooth boundary point. In a sufficiently
small neighborhood of $p$, the nearest-point projection $\pi$ onto $\partial D$
is well defined. For $z$ in this neighborhood, let $n_{\pi(z)}$ be the inward
unit normal and, for $X\in\C^n$, put
\[
 X_z:=\langle X,n_{\pi(z)}\rangle,
 \qquad
 \delta_D(z):=\dist(z,\partial D).
\]
Following \cite{NT2}, set, for $z,w$ close to $p$,
\[
 \widetilde H_D(z,w)
 :=\frac{|(z-w)_z|}
 {\sqrt{|(z-w)_z|}+|z-w|+\sqrt{\delta_D(z)}}+|z-w|.
\]
We shall use the equivalent, slightly simpler scale
\[
 H_D(z,w)
 :=\frac{|(z-w)_z|}
 {\sqrt{|(z-w)_z|}+\sqrt{\delta_D(z)}}+|z-w|.
\]
A simple calculation gives
\[
 \widetilde H_D(z,w)\le H_D(z,w)
 \le \frac43\,\widetilde H_D(z,w).
\]
Thus $H_D\asymp\widetilde H_D$, and the local quasi-symmetry for
$\widetilde H_D$ from \cite{NT2} yields
\[
 H_D(z,w)\asymp H_D(w,z)\lesssim |z-w|^{1/2}.
\]
The three geometric quantities entering $H_D$,
$|(z-w)_z|$, $\delta_D(z)$, and $|z-w|$, reflect genuinely different
boundary regimes and, in general, none of them can be omitted. At a
non-semipositive boundary point, the estimates below show that their
combination gives the sharp local scale. We use $A\lesssim B$ to mean 
$A\le CB$ with a constant uniform in the local configuration under 
consideration, and $A\asymp B$ when both inequalities
hold.

Let $n\ge2$. A point $p\in\partial D$ is called \emph{non-semipositive} if
$\partial D$ is $C^2$-smooth near $p$ and, for a local defining function $\rho$
with $D=\{\rho<0\}$, the Levi form $\LL_\rho(p;\cdot)$ has a negative direction
on $T_p^{\C}\partial D$.

The boundary behaviour of invariant metrics at pseudoconcave and, more
generally, non-semipositive points has been studied in
\cite{Krantz,FL,DNT}. In particular, these works identify the characteristic
normal growth of the infinitesimal Kobayashi, Sibony, and Carath\'eodory
metrics. The present paper is concerned instead with distance-level quantities
and with the pluricomplex Green function. Boundary behaviour of the latter has
been investigated extensively in pseudoconvex settings; see, for example,
\cite{Coman,DH,Herbort}. The non-semipositive situation considered here is of a
different nature, since pseudoconvexity fails locally.

\begin{theorem}\label{thm:two-step}
Let $D\subset\C^n$ be a domain and let $p\in\partial D$ be non-semipositive.
Then, for $z,w$ close to $p$,
\[
 (l_D^*)^{(2)}(z,w)\lesssim H_D(z,w).
\]
\end{theorem}

We now state the corresponding estimates for the pluricomplex Green function.

\begin{theorem}\label{thm:green-lower}
Let $D\subset\C^n$ be a bounded domain and let $p\in\partial D$ be a 
$C^{1,1}$-smooth boundary point. Then, for $z,w$ close to $p$,
\[
 g_D(z,w)\gtrsim H_D(z,w).
\]
\end{theorem}

\begin{theorem}\label{thm:green-upper}
Let $D\subset\C^n$ be a domain and let $p\in\partial D$ be non-semipositive.
Then, for $z,w$ close to $p$,
\[
 g_D(z,w)\lesssim H_D(z,w).
\]
\end{theorem}

Combining Theorem~\ref{thm:two-step}, Theorems~\ref{thm:green-lower} and~
\ref{thm:green-upper}, \cite[Theorems~8 and~9]{NT2}, and \eqref{eq:chain},
gives all the comparisons below.

\begin{corollary}
Let $D\subset\C^n$ be a bounded domain and let $p\in\partial D$ be 
non-semipositive. Then, for $z,w$ close to $p$,
\[
 \begin{aligned}
 g_D(z,w)&\asymp (l_D^*)^{(2)}(z,w)
 \asymp (l_D^{(2)})^*(z,w)\asymp l_D^{(2)}(z,w)\\
 &\asymp k_D^*(z,w)\asymp k_D(z,w)\asymp H_D(z,w).
 \end{aligned}
\]
\end{corollary}

For comparison, if $D$ is bounded and $p$ is non-semipositive, Hartogs extension
across $p$ yields locally
\[
 c_D^*(z,w)\asymp c_D(z,w)\asymp |z-w|;
\]
compare also the corresponding infinitesimal estimates in \cite{DNT}.

The usual Lempert function admits the following one-disc estimate.

\begin{theorem}\label{thm:holder}
Let $D\subset\C^n$, $n\ge2$, and let $p\in\partial D$ be non-semipositive. Then,
for $z,w$ close to $p$,
\[
 l_D^*(z,w)\lesssim\sqrt{H_D(z,w)}\lesssim |z-w|^{1/4}.
\]
\end{theorem}

A matching normal estimate, essentially contained in
\cite[Lemma~6 and Propositions~9 and~12]{NT1}, shows that the exponent $1/4$
is optimal.

\begin{proposition}\label{prop:normal}
Let $D\subset\C^n$ be a bounded domain and let $p\in\partial D$ be a 
$C^{1,1}$-smooth boundary point. Then, for $z,w$ close to $p$ with
$\pi(z)=\pi(w)$,
\[
 l_D^*(z,w)\gtrsim
 \left|\delta_D(z)^{1/4}-\delta_D(w)^{1/4}\right|.
\]
If, moreover, $p$ is non-semipositive, then
\[
 l_D^*(z,w)\asymp
 \left|\delta_D(z)^{1/4}-\delta_D(w)^{1/4}\right|.
\]
\end{proposition}

For comparison, along the same real normal,
\[
 H_D(z,w)\asymp
 \left|\delta_D(z)^{1/2}-\delta_D(w)^{1/2}\right|.
\]

We finally record two qualitative consequences of the same single-disc
geometry.  They show that a point approaching a non-semipositive boundary point
can be joined by a single analytic disc to fixed interior compacta, and that
points approaching two non-semipositive boundary points can be joined with a
uniform parameter bounded away from $1$.

\begin{theorem}\label{thm:qualitative}
Let $D\subset\C^n$, $n\ge2$.

\emph{(i)} Let $p\in\partial D$ be non-semipositive. Then there exists a
neighborhood $U$ of $p$ such that, for every compact $K\Subset D$, there is a
constant $c_K<1$ with
\[
 l_D^*(z,w)\le c_K,
 \qquad z\in K,\quad w\in D\cap U.
\]

\emph{(ii)} Let $p_1,p_2\in\partial D$ be non-semipositive boundary
points. Then there exist neighborhoods $U_j$ of $p_j$ and $c<1$ such that
\[
 l_D^*(z,w)\le c,
 \qquad z\in D\cap U_1,\quad w\in D\cap U_2.
\]
\end{theorem}

Consequently, both conclusions hold also for $g_D$ and for each of the starred
quantities in \eqref{eq:chain}.

\medskip
\noindent\textbf{Outline of the proofs.}
Section~\ref{sec:normal} records the precise
behaviour of the Lempert function along a common normal and shows that the
H\"older exponent $1/4$ is optimal.

Section~\ref{sec:prelim} constructs two families of analytic discs in a
negative Levi direction: fixed-radius discs for the one-disc arguments and
a flexible family for the two-step estimate.

Section~\ref{sec:qualitative} combines the fixed-radius discs with interior
polynomial joining maps. Mergelyan approximation on compact families, with
first-jet interpolation at the boundary-near marked points, preserves the
quadratic interior estimates and gives the two qualitative one-disc bounds.

The estimate $l_D^*\lesssim\sqrt{H_D}$ in Section~\ref{sec:holder} follows
from a two-point Hermite perturbation of a fixed-radius disc.

Section~\ref{sec:two-step} uses flexible discs based at the two endpoints.
Their quadratic corrections are chosen so that they meet at a common
intermediate point at parameter size $O(H_D(z,w))$, yielding the sharp
two-step estimate.

The lower Green estimate in Section~\ref{sec:lower} is independent of
pseudoconcavity. We use a $C^{1,1}$ linearization of the boundary distance
to construct a bounded holomorphic function and glue its logarithmic modulus
into a global plurisubharmonic competitor with the required pole.

The upper Green estimate is first proved on a quadratic model in
Section~\ref{sec:model}, using analytic discs, subharmonicity, and the
Poletsky formula. The uniform interior models of Section~\ref{sec:uniform}
then allow us to transfer this estimate to $D$ in Section~\ref{sec:upper}.

\section{Proof of Proposition~\ref{prop:normal}}\label{sec:normal}

Let $z,w$ be close to $p$ with $\pi(z)=\pi(w)=q$, and, after interchanging
them if necessary, put
\[
 \varepsilon:=\delta_D(z)\ge\delta:=\delta_D(w).
\]
Then
\[
 z=q+\varepsilon n_q,\qquad w=q+\delta n_q.
\]
Choose a boundary patch $\Gamma\Subset\partial D$ containing $p$ and a
fixed tubular neighborhood on which the nearest-point projection is defined.
The proof of \cite[Proposition~12]{NT1} is local and its constants depend only
on the $C^{1,1}$ bounds of the defining functions in normalized charts and on
the size of that tubular neighborhood. These data are uniformly bounded on
$\Gamma$; equivalently, one may cover $\Gamma$ by finitely many of the
neighborhoods used in that proof and take the largest constant. Thus,
uniformly for $q\in\Gamma$,
\[
 l_D^*(z,w)\gtrsim
 \frac{\varepsilon-\delta}{\varepsilon^{3/4}}
 \asymp \varepsilon^{1/4}-\delta^{1/4}.
\]
This proves the first assertion.

Suppose now that $p$ is non-semipositive. Shrink $\Gamma$ so that one
complex-tangential direction has Levi value bounded above by a fixed negative
constant on the whole patch. The inner quadratic model used in the proof of
\cite[Lemma~6]{NT1} can then be chosen with one fixed model size and fixed
quadratic constants: the construction depends only on the uniform $C^2$ norm
of the defining functions and on that negative Levi bound. The corresponding
normal-axis changes of coordinates therefore give a family of
reparametrizations $s_q$ satisfying
\[
 s_q(t)\asymp t,\qquad
 s_q(\varepsilon)-s_q(\delta)\asymp\varepsilon-\delta.
\]
By shrinking the boundary patch if necessary, the constants may be chosen
uniformly in $q$. The model estimate in \cite[Proposition~9]{NT1} has fixed constants.
Hence, by monotonicity and biholomorphic invariance of the Lempert function,
\[
 l_D^*(z,w)
 \lesssim
 \frac{s_q(\varepsilon)-s_q(\delta)}{s_q(\varepsilon)^{3/4}}
 \lesssim
 \frac{\varepsilon-\delta}{\varepsilon^{3/4}},
\]
uniformly for $q\in\Gamma$. Together with the lower estimate, this proves
the second assertion.

\section{Common local preliminaries near a non-semipositive point}\label{sec:prelim}

Fix a $C^2$ defining function $\rho$ near a non-semipositive point $p$, with
$D=\{\rho<0\}$ locally and $|\nabla\rho|=1$ on the boundary. Choose first a unit vector
$v_p\in T_p^{\C}\partial D$ in a strictly negative Levi direction.
If $P_q$ denotes the Hermitian projection onto
$T_q^{\C}\partial D$, then, after shrinking the boundary patch,
$P_qv_p\ne0$ and we may set
\[
 v_q:=\frac{P_qv_p}{|P_qv_p|}.
\]
Since the unit normal is $C^1$ on a $C^2$ boundary, $q\mapsto P_q$ is
Lipschitz on the compact patch; hence so is $q\mapsto v_q$. By continuity
of the Levi form, after shrinking once more there is $\lambda>0$ such that
\begin{equation}\label{eq:neglevi}
 \LL_\rho(q;v_q)\le-4\lambda<0
\end{equation}
throughout the patch.
For $q$ in this patch and $V\in\C^n$, put
\[
 \mathcal S_q(V):=\sum_{j,k}\rho_{z_jz_k}(q)V_jV_k
\]
and define
\[
 A_q^0(V)
 :=-\frac{\mathcal S_q(V)}{2|\partial\rho(q)|^2}\,\overline{\partial\rho(q)}.
\]
Then
\begin{equation}\label{eq:cancel}
 2\sum_j\rho_{z_j}(q)(A_q^0(V))_j+\mathcal S_q(V)=0.
\end{equation}
The map $V\mapsto A_q^0(V)$ is quadratic; on bounded sets it is uniformly
bounded and locally Lipschitz in $V$, with uniformly continuous dependence on
$q$.

\begin{lemma}\label{lem:fixed}
Fix $R>2$. After shrinking the boundary neighborhood and choosing $t>0$
sufficiently small, there are constants $a,b>0$ such that, for every
$w\in D$ sufficiently close to $p$, writing
\[
 w=q+s n_q,
 \qquad q=\pi(w),\qquad s=\delta_D(w),
\]
the map
\[
 F_w(\eta)
 :=q+s n_q+t\eta v_q+t^2\eta^2A_q^0(v_q),
 \qquad |\eta|<R,
\]
satisfies
\begin{equation}\label{eq:fixed-margin}
 \rho(F_w(\eta))\le-as-bt^2|\eta|^2,
 \qquad |\eta|\le R.
\end{equation}
\end{lemma}

\begin{proof}
For the boundary-centered curve
\[
 \Gamma_q(\eta)=q+t\eta v_q+t^2\eta^2A_q^0(v_q),
\]
the linear term in $\rho(\Gamma_q(\eta))$ vanishes because $v_q$ is complex
tangent, while \eqref{eq:cancel} cancels the pure holomorphic quadratic term.
Thus the only quadratic contribution of order $t^2$ in $\eta,\bar\eta$ is
$t^2|\eta|^2\LL_\rho(q;v_q)$. By \eqref{eq:neglevi}, uniformly for $q$ close
to $p$,
\[
 \rho(\Gamma_q(\eta))\le-3\lambda t^2|\eta|^2,
 \qquad |\eta|\le R,
\]
provided $t$ is sufficiently small. Translating by $s n_q$ contributes a
negative term comparable to $-s$. The remaining mixed errors have the form
\[
 Cst|\eta|+Cs\,t^2|\eta|^2+Cs^2.
\]
Using
\[
 st|\eta|\le\varepsilon t^2|\eta|^2+C_\varepsilon s^2
\]
and then shrinking the neighborhood so that $s$ is small, these errors are
absorbed. This proves \eqref{eq:fixed-margin}.
\end{proof}

\begin{lemma}\label{lem:flexible}
Fix $B_0,K_0>0$. There is $M_0=M_0(B_0,K_0)$ such that, for every
$M\ge M_0$, one can first choose $t=t(M)>0$ and then a neighborhood
$U_M$ of $p$ with the following property. Let
\[
 x=q+d\,n_q\in D,
 \qquad d=\delta_D(x),
\]
and let
\[
 U=u+\beta n_q,
 \qquad
 u=Mv_q+h,
 \qquad
 h\in T_q^{\C}\partial D,
 \qquad |h|\le2,
\]
with $|\beta|\le B_0\sqrt d$. If $E\in\C^n$ and $|E|\le K_0M$, then
\[
 F(\eta)
 =x+t\eta U+t^2\eta^2\bigl(A_q^0(U)+E\bigr),
 \qquad |\eta|<1,
\]
maps $\D$ into $D$.
\end{lemma}

\begin{proof}
By \eqref{eq:neglevi} and boundedness of the Levi form,
\[
 \LL_\rho(q;u)\le-cM^2+CM+C.
\]
Choose $M_0$ so large, depending only on $B_0,K_0$ and the fixed local
bounds, that for every $M\ge M_0$, uniformly for all admissible $q,h$,
\[
 \LL_\rho(q;u)\le-\mu_M,
 \qquad \mu_M\ge cM^2.
\]
Taylor's formula at $q$, together with \eqref{eq:cancel}, gives uniformly for
$|\eta|\le1$,
\begin{align*}
 \rho(F(\eta))
 &\le -c_0d+C_0t|\beta||\eta|-c_1\mu_Mt^2|\eta|^2
      +C_1K_0Mt^2|\eta|^2 \\
 &\quad +o(t^2|\eta|^2)
      +O\bigl(dt|\eta|+t^2|\eta|^2|\beta|+d^2\bigr).
\end{align*}
The cancellation \eqref{eq:cancel} removes the pure holomorphic quadratic
term associated with $U$. The mixed Levi terms involving the normal component
of $U$ are $o(1)t^2|\eta|^2$ after the neighborhood is decreased. Moreover,
\[
 C_0t|\beta||\eta|
 \le \frac{c_1\mu_M}{4}t^2|\eta|^2
      +\frac{C_2}{\mu_M}|\beta|^2.
\]
Since $|\beta|\le B_0\sqrt d$, the last term is absorbed by $-c_0d$ for
large $M$, while the $K_0M$ term is absorbed by the negative Levi term.
For this fixed $M$, choose $t=t(M)>0$ sufficiently small and then decrease
the boundary neighborhood to $U_M$ so that all remaining terms are absorbed.
Thus $\rho(F(\eta))<0$ for $|\eta|<1$. This order of choices is the one used
below: $B_0,K_0$ are fixed first, then $M$, then $t$, and only afterwards the
boundary neighborhood.
\end{proof}

\begin{lemma}\label{lem:H}
Let $z,w$ be sufficiently close to $p$, and put
\[
 X:=w-z,\qquad \tau:=H_D(z,w),\qquad
 q_z:=\pi(z),\quad q_w:=\pi(w),
\]
\[
 d_z:=\delta_D(z),\quad d_w:=\delta_D(w),\qquad
 N_z:=\langle X,n_{q_z}\rangle,\quad
 N_w:=\langle X,n_{q_w}\rangle.
\]
Then
\[
 |X|\le\tau,
 \qquad
 |N_z|\lesssim\tau\sqrt{d_z}+\tau^2,
 \qquad
 |N_w|\lesssim\tau\sqrt{d_w}+\tau^2,
\]
and
\[
 |q_z-q_w|+|n_{q_z}-n_{q_w}|\lesssim\tau.
\]
\end{lemma}

\begin{proof}
The first estimate is immediate from the definition of $H_D$. For the second,
write $a=|N_z|$. From the definition of $H_D$,
\[
 a\le \tau\bigl(\sqrt a+\sqrt{d_z}\bigr).
\]
Young's inequality gives $a\lesssim\tau\sqrt{d_z}+\tau^2$. The nearest-point
projection and the normal field are Lipschitz on a tubular neighborhood,
which gives the geometric estimate. Finally, apply the same argument with
$z,w$ interchanged and use the quasi-symmetry $H_D(w,z)\asymp H_D(z,w)$
to obtain the estimate for $N_w$.
\end{proof}

\section{Proof of Theorem~\ref{thm:qualitative}}\label{sec:qualitative}

A \emph{planar chain} is a finite union of compact Jordan domains with
disjoint interiors, consecutive pieces meeting at one node, no other
intersections, and connected complement in $\C$.

\begin{lemma}\label{lem:gluing}
Let $K$ be a fixed planar chain and let $b_1,\ldots,b_m$ ($m\le2$) be
fixed interior points of its terminal pieces. Let $d,d'\in\operatorname{int}K$
be fixed, possibly among the $b_j$. Suppose that
$\overline{\mathcal F}\subset A(K,\C^n)$ is compact, where $A(K,\C^n)$
consists of the continuous maps which are holomorphic on $\operatorname{int}K$.
Assume that $\mathcal F\subset\overline{\mathcal F}$ satisfies
\[
 \rho_j(f(\zeta))\le-a s_j(f)-b|\zeta-b_j|^2,
 \qquad |\zeta-b_j|\le2r,
\]
on fixed disjoint closed discs contained in $\operatorname{int}K$, where
$\rho_j$ are local $C^1$ defining functions, $s_j(f)>0$, and $a,b>0$ are
uniform. Away from any sufficiently small fixed discs about the $b_j$, the
images of $\overline{\mathcal F}$ are relatively compact in $D$.
Then there is $c<1$ such that every $f\in\mathcal F$ admits
$F\in\OO(\D,D)$ with
\[
 F(0)=f(d),\qquad F(\alpha)=f(d'),\qquad |\alpha|\le c.
\]
\end{lemma}

\begin{proof}
Mergelyan approximation can be made uniform on the compact family
$\overline{\mathcal F}$ while preserving the first jets at the $b_j$ and
the values at $d,d'$. Indeed, let $H_f$ be the corresponding interpolation
polynomial and let $Q$ have double zeros at the $b_j$ and simple zeros at
the other marked points. The quotients $(f-H_f)/Q$, with removable values
at the marked points, form a compact subset of $A(K,\C^n)$ by Cauchy's
estimates. A finite cover and Mergelyan's theorem therefore give, for each
fixed $\varepsilon>0$, polynomials $P_f$ of uniformly bounded degree with
\[
 \|P_f-f\|_K<\varepsilon
\]
and the prescribed jets and values.
Since $P_f-f$ has a double zero at $b_j$, Schwarz's lemma on
$\D(b_j,2r)$ gives
\[
 |P_f(\zeta)-f(\zeta)|\le C\varepsilon|\zeta-b_j|^2,
 \qquad |\zeta-b_j|\le r.
\]
Thus, for small $\varepsilon$, the quadratic negative estimates are
preserved. Compactness keeps the rest of $P_f(K)$ inside $D$.

The uniform degree and norm bounds give a uniform derivative bound on a
fixed neighborhood of $K$. The node values range in a compact subset of
$D$ by the hypothesis away from the $b_j$. Hence fixed small node discs
are also mapped into $D$. Inside their union with the interiors of the
pieces, choose a fixed simply connected domain $\Omega$ containing all
marked points. Composing $P_f$ with a Riemann map $\D\to\Omega$ and a
disc automorphism gives the required disc at a fixed parameter
$|\alpha|<1$.
\end{proof}

Fix $r=1/2$ and put $a_w=F_w(r)$, with $F_w$ from
Lemma~\ref{lem:fixed}. By \eqref{eq:fixed-margin}, after shrinking a
neighborhood $U$ of $p$,
\begin{equation}\label{eq:Lcompact}
 L:=\overline{\{a_w:w\in D\cap U\}}\Subset D.
\end{equation}
The neighborhood $U$ depends only on the local boundary geometry.

For compact sets $L_1,L_2\Subset D$, finitely many continuous polynomial
families join every pair in $L_1\times L_2$ through a fixed compact subset
of $D$. To see this, approximate an interior path from $z_0$ to $q_0$ by
a polynomial $P$ with these endpoint values and use
\[
 P_{z,q}(\zeta)=P(\zeta)+(1-\zeta)(z-z_0)+\zeta(q-q_0).
\]
For nearby endpoints this maps a fixed thin neighborhood of $[0,1]$
compactly into $D$; a finite cover of $L_1\times L_2$ gives the assertion.

For part \emph{(i)}, apply this construction to $K_0\times L$, with
$K_0\Subset D$. Attach $F_w$ to the joining map at $a_w=F_w(r)$, using
one of finitely many fixed two-piece planar chains. The parameters
corresponding to $z$ and $w$ lie in the interiors of their pieces. After
a finite closed refinement of the endpoint cover, these piecewise
holomorphic maps have compact closure in $A(K,\C^n)$, including the
limits $w=q+s n_q$ with $s=0$. Near the parameter $b$ corresponding to
$w$, the affine reparametrization of \eqref{eq:fixed-margin} gives
\[
 \rho(f(\zeta))\le-a_1\delta_D(w)-b_1|\zeta-b|^2.
\]
Away from $b$, the images stay in a fixed compact subset of $D$, by the
joining construction and \eqref{eq:fixed-margin}. Lemma~\ref{lem:gluing}
applies on each chain; taking the largest of finitely many parameter
bounds gives $l_D^*(z,w)\le c_{K_0}<1$.

For part \emph{(ii)}, use the fixed-radius discs at $p_1,p_2$ and the
corresponding compact sets $L_1,L_2$. Join them by the interior polynomial
families and attach a local disc at each end. On these finitely many
three-piece chains, \eqref{eq:fixed-margin} gives quadratic negative
estimates at both boundary-near marked points and compact images away
from them. Allowing zero boundary depth gives compact closures of the
families. Lemma~\ref{lem:gluing} yields $l_D^*(z,w)\le c<1$ uniformly.

\section{Proof of Theorem~\ref{thm:holder}}\label{sec:holder}

\begin{proof}
Use the notation of Lemma~\ref{lem:H}, with $\tau=H_D(z,w)>0$.
Fix the radius in Lemma~\ref{lem:fixed} larger than $3$, and put
\[
 U=v_{q_z},\qquad A=A_{q_z}^0(U),\qquad
 f_0(\eta)=z+t\eta U+t^2\eta^2A.
\]
For the fixed $t>0$ supplied by that lemma,
\begin{equation}\label{eq:hermite-base-margin}
 \rho(f_0(\eta))\le-a_0d_z-b_0|\eta|^2,
 \qquad |\eta|\le2.
\end{equation}
Choose a small $c>0$, then large constants $L,K>1$, as specified below,
and set $\alpha=K\sqrt\tau<1/4$. Define
\[
 B(\zeta)=\zeta(\zeta-\alpha),\qquad
 Q_0(\zeta)=\zeta(\zeta-\alpha)^2,\qquad
 Q_1(\zeta)=\zeta^2(\zeta-\alpha),
\]
and
\[
 \gamma_0=
 \begin{cases}
 -N_z/\alpha^3,&d_z<L\tau^2,\\
 0,&d_z\ge L\tau^2,
 \end{cases}
 \qquad
 \gamma_1=
 \begin{cases}
 -\dfrac{t}{\alpha}\langle U,n_{q_w}\rangle
 -\dfrac{N_w}{\alpha^3},&d_w<L\tau^2,\\
 0,&d_w\ge L\tau^2.
 \end{cases}
\]
The required interpolation map is
\begin{equation}\label{eq:hermite-map}
 F(\zeta)=f_0(B(\zeta))+\frac{\zeta}{\alpha}X
 +\gamma_0Q_0(\zeta)n_{q_z}+\gamma_1Q_1(\zeta)n_{q_w}.
\end{equation}
It satisfies $F(0)=z$ and $F(\alpha)=w$. Lemma~\ref{lem:H} gives
$|\langle U,n_{q_w}\rangle|\lesssim\tau$ and
$N_j=O_L(\tau^2)$ whenever $d_j<L\tau^2$. Consequently
\begin{equation}\label{eq:hermite-gamma-est}
 |\gamma_0|+|\gamma_1|\lesssim_L K^{-2}\alpha.
\end{equation}
The chosen corrections cancel the normal first jet at each endpoint
whose boundary depth is less than $L\tau^2$.

\smallskip
\noindent\emph{Endpoint neighborhoods.}
Put $\xi_0=0$, $\xi_1=\alpha$, $x_0=z$, $x_1=w$, and
$d_j=\delta_D(x_j)$. For $|s|\le c\alpha$, expansion of
\eqref{eq:hermite-map} gives
\[
 F(\xi_j+s)=x_j+\alpha sV_j+s^2E_j(s),
\]
where
\[
 V_0=-tU+\frac{X}{\alpha^2}+\gamma_0\alpha n_{q_z},
 \qquad
 V_1=tU+\frac{X}{\alpha^2}+\gamma_1\alpha n_{q_w}.
\]
Their tangential parts differ from $-tv_{q_z}$ and $tv_{q_w}$,
respectively, by $O(K^{-2})+O(\tau)$, so their Levi forms are uniformly
negative for $K$ large. For
$\beta_j=\langle V_j,n_{\pi(x_j)}\rangle$, the definitions and
Lemma~\ref{lem:H} give
\[
 \beta_j=0\quad\hbox{if }d_j<L\tau^2,
 \qquad
 |\beta_j|\lesssim(L^{-1/2}+K^{-2})\sqrt{d_j}
 \quad\hbox{otherwise}.
\]

The term $t^2\alpha^2A$ cancels the principal pure holomorphic quadratic
form by \eqref{eq:cancel}. At $q_w$ the change of this form is $o(1)$,
and the extra quadratic term $t s^2U$ has normal component
$O(\tau)|s|^2=O(K^{-2}\alpha^2|s|^2)$. The Hermite corrections contribute
$O_L(K^{-2}\alpha^2|s|^2)$ by \eqref{eq:hermite-gamma-est}; higher powers
of $s$ are controlled by $|s|\le c\alpha$. Taylor's formula thus gives,
uniformly for $j=0,1$,
\begin{align*}
 \rho(F(\xi_j+s))
 &\le-c_0d_j+C_0\alpha|\beta_j||s|-c_1\alpha^2|s|^2\\
 &\quad +(Cc+C_LK^{-2}+o(1))\alpha^2|s|^2
       +O(d_j\alpha|s|+d_j^2).
\end{align*}
Choose $c$ small, then $L$ large and $K$ large. Young's inequality
absorbs the normal linear term, and shrinking the boundary neighborhood
absorbs the remaining errors. Hence
\begin{equation}\label{eq:hermite-endpoints}
 \rho(F(\xi_j+s))\le-c_2d_j-c_3\alpha^2|s|^2<0,
 \qquad |s|\le c\alpha,\quad j=0,1.
\end{equation}

\smallskip
\noindent\emph{The complement.}
If $|\zeta|,|\zeta-\alpha|\ge c\alpha$, then
\[
 |B|\gtrsim\alpha|\zeta|,\qquad |B|\gtrsim\alpha^2.
\]
Put $F_0=f_0\circ B$ and $R=F-F_0$. From Lemma~\ref{lem:H},
\eqref{eq:hermite-gamma-est}, and $Q_0=B^2/\zeta$, $Q_1=\zeta B$,
\[
 |R|\lesssim_L K^{-2}|B|,
 \qquad
 |\langle R,n_{q_z}\rangle|
 \lesssim K^{-2}\sqrt{d_z}|B|+C_LK^{-2}|B|^2.
\]
For the interpolation term in the normal estimate, use explicitly
\[
 \frac{|N_z||\zeta|}{\alpha}
 \lesssim K^{-2}\sqrt{d_z}|B|+K^{-4}|B|^2;
\]
the two Hermite terms are $O_L(K^{-2}|B|^2)$.
All maps remain in the fixed neighborhood of definition of $\rho$ after
shrinking the boundary neighborhood, since $F-F_0=O(\sqrt\tau)$ uniformly
on $\D$. Its bounded second derivatives give
\[
 \rho(F)\le\rho(F_0)+C|\langle R,n_{q_z}\rangle|
       +C(d_z+|B|)|R|+C|R|^2.
\]
Combining this with \eqref{eq:hermite-base-margin} and Young's inequality
yields, for $K$ sufficiently large,
\[
 \rho(F)\le-\frac{a_0}{2}d_z-(b_0-C_LK^{-2})|B|^2<0.
\]
Together with \eqref{eq:hermite-endpoints}, this proves $F(\D)\subset D$.
Thus $l_D^*(z,w)\le\alpha=K\sqrt{H_D(z,w)}$.
\end{proof}

\section{Proof of Theorem~\ref{thm:two-step}}\label{sec:two-step}

Let $z,w$ be sufficiently close to $p$ and use the notation of
Lemma~\ref{lem:H}. If $\tau=0$, then $z=w$ and there is nothing to prove. It
is enough to treat the case $d_z\ge d_w$. Indeed, if $d_w>d_z$, run the construction
below with $w,z$ in place of $z,w$; since $H_D(w,z)\asymp H_D(z,w)$, reversing the two
constructed discs by automorphisms of $\D$ gives the same two-step chain from
$z$ to $w$. We therefore assume $d_z\ge d_w$. Let $P_q$ denote the Hermitian projection
onto $T_q^{\C}\partial D$ and set
\[
 h_z:=\frac{P_{q_z}X}{2\tau},
 \qquad
 h_w:=-\frac{P_{q_w}X}{2\tau}.
\]
By Lemma~\ref{lem:H}, $|h_z|,|h_w|\le1/2$. Let $M>1$ be a parameter,
to be fixed after the constants $B_0,K_0$ below have been determined, and put
\[
 u_z:=Mv_{q_z}+h_z,
 \qquad
 u_w:=Mv_{q_w}+h_w.
\]
The Lipschitz dependence of $v_q$, $P_q$, and $n_q$, together with
Lemma~\ref{lem:H}, gives
\[
 X-\tau(u_z-u_w)=N_zn_{q_z}+O_M(\tau^2).
\]

If $\tau\le\sqrt{d_z}$, set $\beta_z:=N_z/\tau$; otherwise set
$\beta_z:=0$. In both cases
\[
 |\beta_z|\le B_0\sqrt{d_z}
\]
for a fixed $B_0$, and
\[
 X-\tau\bigl((u_z+\beta_zn_{q_z})-u_w\bigr)=O_M(\tau^2).
\]
Indeed, when $\tau>\sqrt{d_z}$, Lemma~\ref{lem:H} gives
$|N_z|=O(\tau^2)$.

Set
\[
 U_z:=u_z+\beta_zn_{q_z},
 \qquad U_w:=u_w,
\]
\[
 A_z:=A_{q_z}^0(U_z),
 \qquad A_w:=A_{q_w}^0(U_w).
\]
For fixed $M$, continuity of the second derivatives of $\rho$ and the
Lipschitz dependence of $v_q,P_q$ imply, after shrinking the neighborhood,
$|A_z-A_w|\le C_M$. Define
\[
 E_z
 :=\frac{X-\tau(U_z-U_w)-\tau^2(A_z-A_w)}{\tau^2}.
\]
The explicit quadratic formula for $A_q^0$ shows that there is a constant
$K_0>0$, independent of $M$, such that, for each fixed $M$, after shrinking
the neighborhood depending on that $M$,
\[
 |E_z|\le K_0M.
\]
Indeed, the changes caused by $h_z,h_w$ are $O(M)$ with a uniform coefficient,
while the $M^2$ part changes by $o(1)M^2$ as $q_z-q_w\to0$ and can therefore
be made $O(M)$ by shrinking the neighborhood after $M$ is fixed; the
contribution of $\beta_z$ is $o(1)M$. We now choose $M$ sufficiently large
for Lemma~\ref{lem:flexible} with the fixed constants $B_0,K_0$, and then
shrink the neighborhood as just described.

Apply Lemma~\ref{lem:flexible} at $z$ and $w$ to
\[
 F_z(\eta)=z+t\eta U_z+t^2\eta^2(A_z+E_z),
\]
\[
 F_w(\eta)=w+t\eta U_w+t^2\eta^2A_w.
\]
For $a:=\tau/t<1$, the definition of $E_z$ gives
\[
 F_z(a)=F_w(a)=:\xi\in D.
\]
The disc $F_z$ gives $l_D^*(z,\xi)\le a$. Let
\[
 \phi_a(\zeta):=\frac{a-\zeta}{1-a\zeta},
\]
so that $\phi_a$ interchanges $0$ and $a$. Then $F_w\circ\phi_a$ maps
$0$ to $\xi$ and $a$ to $w$, and hence $l_D^*(\xi,w)\le a$. Therefore
\[
 (l_D^*)^{(2)}(z,w)
 \le l_D^*(z,\xi)+l_D^*(\xi,w)
 \le2a
 =\frac2tH_D(z,w).
\]
This proves Theorem~\ref{thm:two-step}.

\section{Proof of Theorem~\ref{thm:green-lower}}\label{sec:lower}

There are $r_0,C_0>0$ such that, for $z,\zeta$ close to $p$,
\[
 \left|\delta_D(\zeta)-\delta_D(z)
 -2\re\langle\zeta-z,\bar\partial\delta_D(z)\rangle\right|
 \le C_0|\zeta-z|^2.
\]
Fix $z,w$ and write
\[
 d=\delta_D(z),\quad T=|z-w|,\quad N=|(z-w)_z|,
\]
and
\[
 L_z(\zeta):=-2\langle\zeta-z,\bar\partial\delta_D(z)\rangle.
\]
Then $|L_z(w)|=N$ and
\[
 \re L_z(\zeta)<d+C_0|\zeta-z|^2.
\]
Choose $A>2$, $C_1>C_0+1$, put $S:=\sqrt N+\sqrt d$, and set
\[
 R=A\max\{S,T\},\qquad a=d+C_1R^2,
\]
with the neighborhood chosen so that $R<r_0$. On
$U_z=D\cap B(z,R)$,
\[
 \re(a-L_z)> (C_1-C_0)R^2>0.
\]
Hence
\[
 s_z(\zeta)=\sqrt{a-L_z(\zeta)},\qquad \re s_z>0,
\]
is holomorphic, and
\[
 h_z(\zeta)=\frac{s_z(\zeta)-\sqrt a}{s_z(\zeta)+\sqrt a}
 =-\frac{L_z(\zeta)}{(\sqrt{a-L_z(\zeta)}+\sqrt a)^2}
\]
satisfies $|h_z|<1$ and $h_z(z)=0$.

Choose $\lambda>0$ with $\lambda\operatorname{diam}D<1$ and
$0<\gamma<\min\{\lambda,1/r_0\}$. Define
\[
 u_z(\zeta)=
 \begin{cases}
 \max\{\log(\gamma R|h_z(\zeta)|),\log(\lambda|\zeta-z|)\},&\zeta\in U_z,\\
 \log(\lambda|\zeta-z|),&\zeta\notin U_z.
 \end{cases}
\]
On $D\cap\partial B(z,R)$ the first entry is at most the second one, so
the psh gluing lemma gives $u_z\in\PSH(D)$. Moreover $u_z<0$. Since
$h_z$ is holomorphic near $z$ and $h_z(z)=0$,
\[
 |h_z(\zeta)|\le C_z|\zeta-z|,\qquad \zeta\to z,
\]
and therefore
\[
 \log(\gamma R|h_z(\zeta)|)
 \le \log|\zeta-z|+O(1).
\]
Thus $u_z(\zeta)\le\log|\zeta-z|+O(1)$ at $z$, so $u_z$ is admissible
for $G_D(z,\cdot)$.

Since $w\in U_z$,
\[
 g_D(z,w)\ge \max\{\gamma R|h_z(w)|,\lambda T\}.
\]
Now $N,d\le R^2/A^2$, hence $a\asymp R^2$ and
\[
 |h_z(w)|\gtrsim \frac{N}{R^2}.
\]
Therefore
\[
 g_D(z,w)\gtrsim \frac NR+T.
\]
If $T\le S$, then $R=AS$ and hence
$N/R+T=A^{-1}N/(\sqrt N+\sqrt d)+T\gtrsim H_D(z,w)$. If $T>S$, then
$N/(\sqrt N+\sqrt d)\le\sqrt N\le S<T$, so $H_D(z,w)<2T$; since the
preceding estimate gives $g_D(z,w)\gtrsim T$, the desired bound follows
in this case as well. This proves Theorem~\ref{thm:green-lower}.

\section{Quadratic model for the Green function}\label{sec:model}

For $\varepsilon>0$ write
\[
 \D_\varepsilon:=\{\lambda\in\C:\ |\lambda|<\varepsilon\}.
\]
Write $\zeta=(\zeta_1,\zeta_2,\zeta')$, where
$\zeta'=(\zeta_3,\ldots,\zeta_n)$ (with the obvious interpretation when
$n=2$). Fix
\[
 \Omega=\{\zeta\in\D_\varepsilon^n:\ \rho(\zeta)<0\},\qquad
 \rho(\zeta)=\re\zeta_1+\kappa|\zeta_1|^2-a|\zeta_2|^2+\mathcal Q(\zeta'),
\]
where $a>0$, $\kappa\ge0$, and $\mathcal Q\ge0$ is a fixed Hermitian form.
The constants below are uniform when $a$ is bounded above and below away from
zero and $\kappa,\mathcal Q$ are bounded above.

\begin{lemma}\label{lem:green-model}
For
\[
 P=(-d,0,0),\qquad W=(-d+s,\eta,\xi)\in\Omega,
\]
put $N=|s|$ and $T=|W-P|$. For $P,W$ sufficiently close to the origin,
\begin{equation}\label{eq:model-est}
 g_\Omega(P,W)\lesssim
 \frac{N}{\sqrt N+\sqrt d}+T.
\end{equation}
\end{lemma}

\begin{proof}
Fix a sufficiently small constant $\tau_0>0$. We separate two regimes.

\smallskip\noindent
\emph{I. The regime $d\le \tau_0^{-1}(\sqrt N+T)^2$.}
We first prove
\[
 g_\Omega(P,W)\lesssim \sqrt N+T
\]
under the above assumption, with constants allowed to depend on the fixed
number $\tau_0$. In this regime $\sqrt d\lesssim\sqrt N+T$, and hence the
right hand side of \eqref{eq:model-est} is comparable with $\sqrt N+T$.
Indeed, writing $r=\sqrt N+T$, one has
\[
 \frac{N}{\sqrt N+\sqrt d}+T
 \gtrsim \frac{N}{r}+T\asymp r,
\]
while the reverse inequality is immediate.

Put $r=\sqrt N+T$ and $R=Ar$, where $A\gg1$ is fixed.  Fix
$b=a/4$.  Write $w_1=-d+s$.  After shrinking the model once and for all,
\[
 |1+2\kappa\overline{w_1}|\ge \frac12.
\]
Define
\[
 A_W:=\frac{2a\bar\eta R}{1+2\kappa\overline{w_1}},
 \qquad
 B_W:=\frac{bR^2}{1+2\kappa\overline{w_1}},
\]
and
\[
 \begin{aligned}
 \Phi_1(\lambda)&=w_1+A_W\lambda+B_W\lambda^2,\\
 \Phi_2(\lambda)&=\eta+R\lambda,\\
 \Phi_{\mathrm{tail}}(\lambda)&=\xi.
 \end{aligned}
\]
Then $\Phi(0)=W$.  Put $\Delta(\lambda)=A_W\lambda+B_W\lambda^2$.
Using
\[
 \re(w_1+\Delta)+\kappa|w_1+\Delta|^2
 =\re w_1+\kappa|w_1|^2
  +\re\!\left((1+2\kappa\overline{w_1})\Delta\right)
  +\kappa|\Delta|^2,
\]
the definitions of $A_W,B_W$ give the exact cancellation
\[
 \re\!\left((1+2\kappa\overline{w_1})A_W\lambda\right)
 =2a\re(\bar\eta R\lambda)
\]
and
\[
 \re\!\left((1+2\kappa\overline{w_1})B_W\lambda^2\right)
 =bR^2\re\lambda^2.
\]
Consequently
\begin{align*}
 \rho(\Phi(\lambda))
 &=\rho(W)+bR^2\re\lambda^2-aR^2|\lambda|^2
   +\kappa|\Delta(\lambda)|^2\\
 &\le \rho(W)-(a-b)R^2|\lambda|^2+C R^4|\lambda|^2.
\end{align*}
Here we used $|A_W|\lesssim|\eta|R\lesssim R^2/A$ and
$|B_W|\lesssim R^2$.  Shrinking the fixed local scale so that
$CR^2\le(a-b)/2$, we obtain the uniform estimate
\[
 \rho(\Phi(\lambda))
 \le \rho(W)-\frac{a-b}{2}R^2|\lambda|^2<0,
 \qquad \lambda\in\D.
\]
Since $|\Phi_2(\lambda)|\le |\eta|+R$, $|\Phi_{\mathrm{tail}}(\lambda)|=|\xi|$,
and $|\Phi_1(\lambda)|\le |w_1|+C R^2$, the truncating polydisc is also
respected after shrinking the fixed local scale. Thus $\Phi(\D)\subset\Omega$.

For $|\lambda|=1$ write
\[
 W_\theta=\Phi(e^{i\theta})=(-d+\sigma_\theta,\eta_\theta,\xi),
\]
where
\[
 \eta_\theta=\eta+Re^{i\theta},\qquad
 \sigma_\theta=s+A_We^{i\theta}+B_We^{2i\theta}.
\]
Since $|\eta|,|\xi|\le T\le R/A$ and $N\le R^2/A^2$, while
$|A_W|\lesssim R^2/A$ and, after shrinking the local scale,
$|B_W|\le(a/2)R^2$. Indeed,
\[
 a|\eta_\theta|^2-\mathcal Q(\xi)
 \ge aR^2-2aR|\eta|-a|\eta|^2-\mathcal Q(\xi)
 \ge aR^2-\frac{C}{A}R^2,
\]
whereas
\[
 |\sigma_\theta|
 \le |s|+|A_W|+|B_W|
 \le \frac{R^2}{A^2}+\frac{C}{A}R^2+\frac a2R^2.
\]
Thus, after choosing $A$ sufficiently large, there is a uniform
$c>0$ such that
\begin{equation}\label{eq:margin}
 |\sigma_\theta|
 \le a|\eta_\theta|^2-\mathcal Q(\xi)-cR^2
\end{equation}
uniformly in $\theta$.

Put $\alpha=LR$ with $L\gg1$ fixed. After shrinking the local scale once more,
we may assume $\alpha<1$. Define
\[
 \Psi_\theta(\zeta)=
 \left(-d+\frac{\sigma_\theta}{\alpha^2}\zeta^2,
 \frac{\eta_\theta}{\alpha}\zeta,
 \frac{\xi}{\alpha}\zeta\right).
\]
Then $\Psi_\theta(0)=P$ and $\Psi_\theta(\alpha)=W_\theta$.
By \eqref{eq:margin},
\[
 \re\!\left(\frac{\sigma_\theta}{\alpha^2}\zeta^2\right)
 -a\left|\frac{\eta_\theta}{\alpha}\zeta\right|^2
 +\mathcal Q\!\left(\frac{\xi}{\alpha}\zeta\right)
 \le -\frac{c}{L^2}|\zeta|^2.
\]
Moreover, since $d=O(R^2)$ and $|\sigma_\theta|=O(R^2)$,
\[
 \kappa|\Psi_{\theta,1}(\zeta)|^2
 \le C d^2+C\frac{d}{L^2}|\zeta|^2
      +C\frac1{L^4}|\zeta|^4.
\]
Therefore
\[
 \rho(\Psi_\theta(\zeta))
 \le -d+Cd^2
 -\left(\frac{c}{L^2}-C\frac d{L^2}-\frac C{L^4}\right)|\zeta|^2.
\]
First choose $L$ so large that $C/L^4\le c/(4L^2)$, and then shrink
the local scale so that $Cd\le\min\{1/2,c/4\}$.  It follows that
\[
 \rho(\Psi_\theta(\zeta))
 \le -\frac d2-\frac{c}{2L^2}|\zeta|^2<0,
 \qquad |\zeta|<1.
\]
For $L$ large the truncating polydisc is also respected. Holomorphic
contractibility gives
\[
 g_\Omega(P,W_\theta)\le |\alpha|\lesssim R.
\]
Since $G_\Omega(P,\Phi(\lambda))$ is subharmonic,
\[
 G_\Omega(P,W)\le \frac1{2\pi}\int_0^{2\pi}
 G_\Omega(P,W_\theta)\,d\theta\le\log(CR),
\]
and hence $g_\Omega(P,W)\lesssim R\asymp\sqrt N+T$.

\smallskip\noindent
\emph{II. The remaining regime.}
It remains to consider
\[
 (\sqrt N+T)^2<\tau_0 d.
\]
In particular,
\begin{equation}\label{eq:smallregime}
 N\le \tau_0 d,\qquad T\le \sqrt{\tau_0}\,\sqrt d.
\end{equation}
We choose the fixed constant $\tau_0>0$ sufficiently small in the estimates
below.

If $W=P$, there is nothing to prove, so assume $W\ne P$.

Set
\[
 \beta=K\left(\frac N{\sqrt d}+T\right),\qquad A_0=\sqrt\beta,
\]
where $K>1$ will be chosen sufficiently large. If $\beta\ge\beta_0$ for
some fixed $\beta_0>0$, then $g_\Omega\le1\le\beta_0^{-1}\beta$, so there
is nothing to prove. Thus we may assume $\beta\ll1$. For an integer
$m\ge1$, to be chosen below, let
\[
 \chi(\lambda)=\frac{A_0-\lambda^m}{1-A_0\lambda^m}.
\]
Since $0<A_0<1$, the function $\chi$ is a holomorphic self-map of
$\D$, in particular $|\chi|<1$ on $\D$. Moreover $\chi(0)^2=\beta$, and its
$m$ zeros $\lambda_j$ satisfy $|\lambda_j|=\beta^{1/(2m)}$. Define
\[
 \begin{aligned}
 \Theta_1(\lambda)&=-d+\frac{s}{\beta}\chi(\lambda)^2,\\
 \Theta_2(\lambda)&=\chi(\lambda)^2\left(\frac\eta\beta+\gamma\lambda\right),\\
 \Theta_{\mathrm{tail}}(\lambda)&=\chi(\lambda)^2\frac\xi\beta,
 \end{aligned}
\]
where $0<\gamma<\varepsilon/4$ is fixed. Then $\Theta(0)=W$ and
$\Theta(\lambda_j)=P$ for all $j$.

We check that $\Theta(\D)\subset\Omega$. Put $y=|\chi(\lambda)|^2$. Then
\begin{equation}\label{eq:rhoTheta}
 \rho(\Theta(\lambda))\le
 -d+\frac N\beta y
 -ay^2\left|\frac\eta\beta+\gamma\lambda\right|^2
 +y^2\frac{\mathcal Q(\xi)}{\beta^2}
 +\kappa\left|-d+\frac{s}{\beta}\chi(\lambda)^2\right|^2.
\end{equation}
For $|\lambda|\ge1/2$, the definition of $\beta$, together with
$|\eta|,|\xi|\le T$, gives
\[
 \frac{|\eta|}{\beta}\le\frac1K,\qquad
 \frac{\mathcal Q(\xi)}{\beta^2}\lesssim\frac1{K^2},\qquad
 \frac N\beta\le\frac{\sqrt d}{K}.
\]
Choose first $K$ large enough that $|\eta|/\beta\le\gamma/4$ and
$\mathcal Q(\xi)/\beta^2$ is smaller than half of the resulting negative
quadratic coefficient. Then
\[
 \rho(\Theta(\lambda))\le-d+\frac N\beta y-cy^2
 \le-d+C'\left(\frac N\beta\right)^2<0,
\]
because $\beta\ge K N/\sqrt d$.  The last term in
\eqref{eq:rhoTheta} is bounded by $C(d^2+d/K^2)$ and is absorbed after
increasing $K$ and shrinking the local scale.
For $|\lambda|\le1/2$, fix a small $\delta>0$ and choose
$m=m(\beta,\delta)$ sufficiently large
that $2^{-m}\ll_\delta\sqrt\beta$; then
$|\chi(\lambda)|^2\le(1+\delta)\beta$. Dropping the negative tangential term in
\eqref{eq:rhoTheta} and using $|\chi|^2\le(1+\delta)\beta$, the
$\kappa$-term is $O(d^2+N^2)$. Hence
\[
 \rho(\Theta(\lambda))
 \le-d+(1+\delta)N+(1+\delta)^2\mathcal Q(\xi)+O(d^2+N^2)<0
\]
by \eqref{eq:smallregime}, after choosing $\tau_0,\delta$ sufficiently small
and shrinking the local scale. The truncating polydisc follows from
\[
 |\Theta_1+d|\le N/\beta,\qquad
 |\Theta_2|\le |\eta|/\beta+\gamma,\qquad
 |\Theta_{\mathrm{tail}}|\le |\xi|/\beta.
\]
Thus $\Theta(\D)\subset\Omega$.

By the disc formula for the additive pluricomplex Green function, after
exponentiation (see \cite[Theorem~1]{EdiGreen}; cf. also
\cite[Theorem~17.4.3]{JP}),
\[
 g_\Omega(P,W)
 \le
 \prod_{\lambda\in\Theta^{-1}(P)}
 |\lambda|^{\operatorname{ord}_{\lambda}(\Theta-P)}.
\]
At each of the $m$ zeros $\lambda_j$ of $\chi$, all components of
$\Theta-P$ contain the factor $\chi^2$, hence
$\operatorname{ord}_{\lambda_j}(\Theta-P)\ge2$.  Any additional
preimages or higher multiplicities can only decrease the product.
Consequently
\[
 g_\Omega(P,W)\le\prod_{j=1}^m|\lambda_j|^2=\beta
 \lesssim \frac N{\sqrt d}+T.
\]
Under \eqref{eq:smallregime} the denominator in \eqref{eq:model-est} is
comparable with $\sqrt d$, proving the lemma.
\end{proof}

\section{Uniform quadratic models for the Green function}\label{sec:uniform}

\begin{lemma}\label{lem:uniform-model}
Let $p\in\partial D$ be a non-semipositive boundary point.
After shrinking a neighborhood $U$ of $p$, there are constants
$\varepsilon,a,\kappa,\Lambda,t_0>0$ such that for every $q\in\partial D\cap U$
there is a local biholomorphic coordinate map $\zeta=\Phi_q(Z)$, defined on a
neighborhood $V_q$ of uniform size, with $\Phi_q(q)=0$. Its differential
$d\Phi_q(q)$ is unitary and sends the complex normal direction to the first
coordinate direction. If $n_q$ denotes the inward unit normal, the inward normal
segment is normalized by
\[
 \Phi_q(q+t n_q)=(-t,0,0),\qquad 0<t<t_0.
\]
The maps and their inverses have uniformly bounded $C^2$ norms and may be
chosen to depend continuously on $q$ in the $C^2$ topology on these uniform
coordinate neighborhoods. In these coordinates
\begin{equation}\label{eq:uniformmodel2}
 \Omega_q:=\left\{\zeta\in\D_\varepsilon^n:
       \re\zeta_1+\kappa|\zeta_1|^2
       -a|\zeta_2|^2+\Lambda|\zeta'|^2<0\right\}
 \subset \Phi_q(D\cap V_q).
\end{equation}
The same constants work for every $q\in\partial D\cap U$.
\end{lemma}

\begin{proof}
Use the defining function $\rho$ and the negative Levi field $v_q$ fixed
above. After shrinking $U$ if necessary,
\[
 \LL_\rho(q;v_q)\le-4\mu<0,
 \qquad q\in\partial D\cap U,
\]
for some $\mu>0$. Complete $v_q$ and the complex normal $e_1(q)$ to a
unitary frame.
Taylor expansion in the corresponding affine unitary coordinates gives,
uniformly in $q$,
\[
 \rho_q(\zeta)=\re\zeta_1+\mathcal R_q(\zeta)+o(|\zeta|^2),
\]
where $\mathcal R_q$ is a real quadratic polynomial.  We now make a holomorphic
quadratic shear in the first coordinate, involving only tangential
variables.  This shear is tangent to the identity and hence preserves the
normal axis $\{(\zeta_1,0):\zeta_1\in\C\}$.  It removes the pure
holomorphic quadratic terms involving only tangential variables.  The
remaining quadratic terms which contain $\zeta_1$ are bounded by
\[
 K|\zeta_1|^2+\mu|\zeta_2|^2+C_1|\zeta'|^2,
\]
with uniform $K,C_1$.  Let
\[
 \mathcal H_q(\zeta_T):=\sum_{j,k=2}^n
 \frac{\partial^2 \rho_q}{\partial \zeta_j\,\partial\bar\zeta_k}(0)
 \zeta_j\bar\zeta_k
\]
be the Hermitian tangential quadratic part.  In the chosen tangential frame,
the strict negative bound in the $v_q$ direction gives
\[
 \mathcal H_q(\zeta_T)
 \le-3\mu|\zeta_2|^2+C_2|\zeta'|^2.
\]
Because the Taylor remainder is $o(|\zeta|^2)$ uniformly in $q$, after shrinking
a fixed polydisc we retain a fixed part of the negative Levi margin and obtain
\begin{equation}\label{eq:rho-upper}
 \rho_q(\zeta)
 \le \re\zeta_1+\kappa|\zeta_1|^2
       -\frac32\mu|\zeta_2|^2+C_3|\zeta'|^2
\end{equation}
with uniform $\kappa,C_3>0$.  Put $a=\mu$ and choose $\Lambda>C_3$.  If
\[
 M_q(\zeta):=
 \re\zeta_1+\kappa|\zeta_1|^2-a|\zeta_2|^2+\Lambda|\zeta'|^2<0,
\]
then \eqref{eq:rho-upper} gives
\[
 \rho_q(\zeta)
 \le M_q(\zeta)-\frac\mu2|\zeta_2|^2-(\Lambda-C_3)|\zeta'|^2<0.
\]
This proves \eqref{eq:uniformmodel2}.  All constants are uniform because
the boundary patch has compact closure, the negative Levi bound is uniform,
and the $C^2$ Taylor remainders are uniform there. The affine unitary frames
and the coefficients of the quadratic shears can be chosen continuously in
$q$; hence, after shrinking the common coordinate size once more, the maps
$\Phi_q$ and their inverses depend continuously on $q$ in the $C^2$ topology
and have uniform $C^2$ bounds. Since the shear involves only tangential
variables, it fixes the normal axis, which gives
$\Phi_q(q+t n_q)=(-t,0,0)$.
\end{proof}

\begin{lemma}\label{lem:model-transition}
In the notation of Lemma~\ref{lem:uniform-model}, write
$\widehat M_q:=M_q\circ\Phi_q$ in the original coordinates. There are a
constant $C>0$ and a modulus of continuity $\omega$, with $\omega(t)\to0$ as
$t\to0$, such that for $q,q'$ in the boundary patch and $Z$ in the fixed
common coordinate neighborhood, with
\[
 h:=|q-q'|,\qquad R_0:=|Z-q'|,
\]
one has
\begin{equation}\label{eq:model-transition}
 \widehat M_q(Z)-\widehat M_{q'}(Z)
 \le C(hR_0+h^2)+\omega(h)R_0^2.
\end{equation}
\end{lemma}

\begin{proof}
Expand both functions at $q'$. The first-order normalization gives
$d\widehat M_q(q)=d\rho(q)$ and
$d\widehat M_{q'}(q')=d\rho(q')$. Since $q,q'\in\partial D$, Taylor's
formula for $\rho$ at $q$ yields
$d\rho(q)(q'-q)=O(h^2)$. The uniform $C^2$ bounds for $\widehat M_q$ and the
$C^1$ regularity of $d\rho$ therefore imply
\[
 \widehat M_{q'}(q')=0,
 \qquad
 \widehat M_q(q')=O(h^2),
\]
and
\[
 |d\widehat M_q(q')-d\widehat M_{q'}(q')|
 \le |d\widehat M_q(q')-d\widehat M_q(q)|
     +|d\rho(q)-d\rho(q')|
 \le Ch.
\]
By the $C^2$-continuous dependence of the normalizations on the base point,
there is a modulus of continuity $\omega$ such that
\[
 \|d^2\widehat M_q-d^2\widehat M_{q'}\|\le\omega(h)
\]
on the fixed common coordinate neighborhood. Subtracting the two Taylor
expansions at $q'$ gives \eqref{eq:model-transition}.
\end{proof}

\section{Proof of Theorem~\ref{thm:green-upper}}\label{sec:upper}

We follow the geometric decomposition used in \cite{NT2} for the upper
estimate for $k_D$, keeping $z$ as the pole throughout.
By Lemma~\ref{lem:uniform-model}, after shrinking $U$,
every $q\in\partial D\cap U$ admits a uniform interior quadratic model
$\Omega_q$ of the form \eqref{eq:uniformmodel2}.  Let $\Phi_q$ denote the
corresponding local biholomorphic coordinate map and put
\[
 \widetilde\Omega_q:=\Phi_q^{-1}(\Omega_q)\subset D.
\]
Thus, whenever $z,w\in\widetilde\Omega_q$, monotonicity and
biholomorphic invariance give
\[
 g_D(z,w)\le g_{\widetilde\Omega_q}(z,w)
 =g_{\Omega_q}\bigl(\Phi_q(z),\Phi_q(w)\bigr).
\]
We use the notation $\pi$, $\delta_D$, and $n_q$ fixed above. For
$u\in D\cap U$,
\[
 u=\pi(u)+\delta_D(u)n_{\pi(u)}.
\]

For $z,w\in D\cap U$ put
\[
 d=\delta_D(z),\quad X=w-z,
\]
and
\[
 Q=Q_D(z,w):=|X_z|+|X|^2+|X|\sqrt d.
\]
The $C^2$ Taylor estimate for the signed distance gives, without
interchanging the two points,
\begin{equation}\label{eq:delta-diff}
 |\delta_D(z)-\delta_D(w)|
 \lesssim |X_z|+|X|^2\lesssim Q.
\end{equation}

\subsection*{Small-Q regime}
Assume $Q\le c_*d$, where $c_*>0$ is sufficiently small. Then
\[
 |X|\le Q/\sqrt d\le c_*\sqrt d,
 \qquad |X_z|\le Q.
\]
Take the coordinates at $q=\pi(z)$ and set
\[
 P:=\Phi_q(z)=(-d,0,0),\qquad W:=\Phi_q(w)=P+Y.
\]
Since the nonlinear shear in the first coordinate contains only
tangential quadratic terms, while the remaining coordinates are unitary to
first order, the normalization gives
\[
 |Y_1|\lesssim |X_z|+|X|^2,\qquad |Y_{\rm tan}|\lesssim|X|.
\]
Consequently the defining function of the model satisfies
\[
 M_q(W)\le-d+C d^2+C Q.
\]
After shrinking $U$, we may assume $Cd\le1/4$; choosing $c_*>0$ still
smaller so that $Cc_*\le1/4$, the assumption $Q\le c_*d$ yields
\[
 M_q(W)\le-\frac12d<0.
\]
Thus $P,W\in\Omega_q$, equivalently $z,w\in\widetilde\Omega_q$, and by
monotonicity, biholomorphic invariance, and Lemma~\ref{lem:green-model},
\[
 g_D(z,w)\le g_{\Omega_q}(P,W)
 \lesssim \frac{|X_z|+|X|^2}{\sqrt d}+|X|
 \lesssim \frac{|X_z|}{\sqrt d}+|X|.
\]
Since $Q\le c_*d$, we have
$\sqrt{|X_z|}+\sqrt d\asymp\sqrt d$; hence the last expression is
comparable with $H_D(z,w)$.

\subsection*{Large-Q regime}
Assume $Q\ge c_*d$ and put $r=\sqrt Q$.  By
\eqref{eq:delta-diff}, both $d$ and $d_w:=\delta_D(w)$ are $O(r^2)$,
while $|z-w|\le r$. Put $q_z=\pi(z)$ and $q_w=\pi(w)$. Then
$|q_z-q_w|\lesssim r$.

We shall use Lemma~\ref{lem:model-transition} for the two boundary anchors
$q_z$ and $q_w$.

Fix now a sufficiently large constant $L$ and put $R=Lr$.  In the
coordinates at $q_w$ consider the analytic disc
\[
 \varphi(\lambda)=\Phi_{q_w}^{-1}(-d_w,R\lambda,0),
 \qquad \lambda\in\D.
\]
It is centered at $w$. Since
\[
 M_{q_w}(-d_w,R\lambda,0)
 =-d_w+\kappa d_w^2-aR^2|\lambda|^2,
\]
after shrinking $U$ (with $L$ fixed) the whole disc lies in
$\widetilde\Omega_{q_w}\subset D$ and the truncating polydisc is respected.
For $W_\theta=\varphi(e^{i\theta})$ we moreover have
\[
 \widehat M_{q_w}(W_\theta)\le-aR^2.
\]

Apply \eqref{eq:model-transition} with $q=q_z$, $q'=q_w$ and
$Z=W_\theta$.  Here $h\lesssim r$ and $|W_\theta-q_w|\lesssim R$.
Hence
\[
 \widehat M_{q_z}(W_\theta)
 \le -aL^2r^2+C(L+1)r^2+\omega(Cr)L^2r^2.
\]
Choose first $L$ so large that $C(L+1)\le aL^2/4$, and then shrink $U$
so that $\omega(Cr)\le a/4$.  Thus
\[
 \widehat M_{q_z}(W_\theta)\le-\frac a2L^2r^2<0.
\]
Since $|W_\theta-q_z|\lesssim R+r$ and $L$ is fixed, the truncating
polydisc for the $q_z$-model is also respected after shrinking $U$ once more.
Hence $W_\theta\in\widetilde\Omega_{q_z}$ for every $\theta$.  Of course
$z\in\widetilde\Omega_{q_z}$ as well.

We also need the model size of $W_\theta-z$.  In the original coordinates,
by the $C^2$ normalization at $q_w$,
\[
 W_\theta-w=R\,v_\theta+O(R^2),
\]
where $v_\theta$ is a unit complex-tangential vector at $q_w$.  Since the
boundary normal is Lipschitz and $|q_z-q_w|\lesssim r$,
\[
 |\langle v_\theta,n_{q_z}\rangle|\lesssim r.
\]
Together with $|X_z|\le Q=r^2$, $|X|^2\le Q$, and the $C^2$
comparison of normal components, this gives, in the coordinates at $q_z$,
\[
 \widetilde T_\theta
 :=|\Phi_{q_z}(W_\theta)-\Phi_{q_z}(z)|\lesssim r,
 \qquad
 \widetilde N_\theta
 :=|\Phi_{q_z}(W_\theta)_1-\Phi_{q_z}(z)_1|\lesssim r^2.
\]
Indeed the normal contribution of the disc displacement is
$O(R|q_z-q_w|+R^2)$, while the contribution of $w-z$ is
$O(|X_z|+|X|^2)$.
Since also $d\lesssim r^2$, Lemma~\ref{lem:green-model}, monotonicity
and biholomorphic invariance yield
\[
 g_D(z,W_\theta)
 \le g_{\Omega_{q_z}}\bigl(\Phi_{q_z}(z),\Phi_{q_z}(W_\theta)\bigr)
 \lesssim r.
\]
Finally $G_D(z,\varphi(\lambda))$ is subharmonic, and therefore
\[
 G_D(z,w)
 \le\frac1{2\pi}\int_0^{2\pi}G_D(z,W_\theta)\,d\theta
 \le\log(C_Lr).
\]
Thus
\[
 g_D(z,w)\lesssim\sqrt{Q_D(z,w)}.
\]

It remains to compare $\sqrt Q$ directly with $H_D(z,w)$. Put
\[
 A:=|X_z|,\qquad t:=|X|,
\]
so that
\[
 Q=A+t^2+t\sqrt d,\qquad d\lesssim Q.
\]
Fix a sufficiently small $\varepsilon>0$. If
$t\ge\varepsilon\sqrt Q$, then
\[
 H_D(z,w)\ge t\gtrsim\sqrt Q.
\]
If $t<\varepsilon\sqrt Q$, then
\[
 t^2+t\sqrt d
 \le(\varepsilon^2+C\varepsilon)Q.
\]
Choosing $\varepsilon$ sufficiently small gives $A\gtrsim Q$.
Moreover $\sqrt d\lesssim\sqrt Q$, and hence
\[
 \frac{A}{\sqrt A+\sqrt d}
 \gtrsim\frac{Q}{\sqrt Q}
 =\sqrt Q.
\]
Thus in all cases
\[
 \sqrt Q\lesssim H_D(z,w).
\]
Therefore $g_D(z,w)\lesssim H_D(z,w)$ also in the large-$Q$ regime.
This proves Theorem~\ref{thm:green-upper}.

\medskip
\noindent\textbf{AI statement.}
ChatGPT (OpenAI) was used in developing parts of the proofs and editing the manuscript.
The authors take full responsibility for the mathematical content of the paper.

\end{document}